\documentclass[12pt]{article}
\usepackage{amssymb}
\usepackage{amsmath}
\usepackage{amsfonts}

\begin{document}

\renewcommand{\theequation}{\arabic{section}.\arabic{equation}}

\begin{center}
{\Large Existence for a free boundary problem describing a propagating
disturbance }

\medskip 

{\Large Gabriela Marinoschi}

\medskip

\textquotedblleft Gheorghe Mihoc-Caius Iacob\textquotedblright\ Institute of
Mathematical Statistics and

Applied Mathematics of the Romanian Academy,

Calea 13 Septembrie 13, Bucharest, Romania

gabriela.marinoschi@acad.ro

\medskip
\end{center}

\noindent \textsc{Abstract.} We prove the existence of a solution to an 1-D
free boundary problem which describes the propagation of disturbances of
shock type, modeled by a non standard variational inequality.

\bigskip


\bigskip


\section{Introduction}

\setcounter{equation}{0}

This paper concerns the existence of a solution to the differential inclusion%
\begin{eqnarray}
\frac{dL}{dt}(t)+\partial I_{K}(\Gamma (t,L(t))) &\ni &U(t,L(t)),\mbox{ a.e. 
}t\in (0,T),  \label{L1} \\
L(0) &=&L_{0}\geq 0,  \notag
\end{eqnarray}%
where $\Gamma $ and $U$ depending on $t$ and $x$ are given functions, $K$ is
the set $\{z\in \mathbb{R};$ $z\geq \Gamma ^{\ast }>0\}$ and $\partial I_{K}$
is the subdifferential of $I_{K},$ the indicator function of $K.$ We note
that (\ref{L1}) can be equivalently written as 
\begin{equation}
\frac{dL}{dt}(t)=U(t,L(t)),\mbox{ \ in }\{t\geq 0;\mbox{ }\Gamma
(t,L(t))>\Gamma ^{\ast }\},  \label{L1-1}
\end{equation}%
\begin{equation}
\frac{dL}{dt}(t)\geq U(t,L(t)),\mbox{ \ in }\{t\geq 0;\mbox{ }\Gamma
(t,L(t))=\Gamma ^{\ast }\}.  \label{L1-1-0}
\end{equation}%
The variational inequality (\ref{L1}) can describe a discontinuity occurring
at the surface $L(t)$ of a system of particles moving with the velocity $%
U(t,x),$ when an intrinsic constraint forces the particles lying on the
surface $L(t)$ to advance with a velocity greater than $U(t,L(t)).$ More
precisely, a particle on the surface $x=L(t)$ moves at each time with the
surface velocity $U(t,L(t))$ as long as the function $\Gamma (t,L(t))$
exceeds a prescribed value $\Gamma ^{\ast },$ but it is pushed out from the
surface with a velocity greater than $U(t,L(t))$ if $\Gamma (t,L(t))$
decreases up to $\Gamma ^{\ast },$ or below it, at a moment $t.$ In other
words, $L(t)$ remains a material surface advancing with the group velocity
as well as $\Gamma (t,L(t))$ is larger than $\Gamma ^{\ast }$ and exhibits a
discontinuous behavior if $\Gamma (t,L(t))$ is equal or lower than $\Gamma
^{\ast }.$

The starting point for the study of such a variational inequality was an 1-D
model of epidermis cell growth introduced in \cite{GIM1} for the stationary
case, and developed in \cite{GIM2} for the dynamical case. We stress that
the current approach (\ref{L1}) of the behavior of the 1-D epidermis free
boundary is different by that treated in \cite{GIM1} and \cite{GIM2}, where
another process modeled by a reversed inequality was studied.

Since the study of (\ref{L1}) is related to this 1-D model, the analysis was
restricted here to the 1-D case. In few words, the growth of the epidermis,
viewed as a body of different type of cells, takes place with the group
velocity $U(t,x)$. The bonds between cells maintain the tissue firm if the
cohesion between cells represented by the function $\Gamma (t,x)$ has values
greater than the critical threshold $\Gamma ^{\ast }.$ Thus, at each time $%
t, $ the cells found at the position $L(t),$ representing the free surface,
advance with the velocity $U(t,L(t)),$ according to (\ref{L1-1}). In this
case, the boundary is material and this assumption was often used in the
literature for a cell system advance (see e.g. \cite{Friedman-1}). By cell
ageing, the cohesion decreases and it is lost when $\Gamma $ reaches $\Gamma
^{\ast }$ at $L(t),$ producing a shock-type action for the surface velocity,
which is a process described by our current model. The particles lying on
the surface $L(t)$ at that moment are detached from the system and are
thrown outside it in the direction of movement, the surface velocity
exhibiting a jump described by (\ref{L1-1-0}).

Another interpretation of the variational inequality (\ref{L1}) may be
related to hysteresis processes, because (\ref{L1}) is similar to the
differential inclusions arising in mathematical models of hysteresis of stop
type used in the study of rheological models (see e.g., \cite{Visintin}, p.
25).

We note that (\ref{L1-1})-(\ref{L1-1-0}) is a free boundary problem, because
the set $\{t;$ $\Gamma (t,L(t))=\Gamma ^{\ast }\}$ and the domain $\{t;$ $%
\Gamma (t,L(t))>\Gamma ^{\ast }\}$ are unknown. We shall prove below the
existence of a solution to (\ref{L1}) in a certain generalized
(distributional) sense. To this purpose, first we introduce some notation
and recall a few definitions. We shall use results from the monographs \cite%
{HB-1} and \cite{HB-2}, but as they already belong to a classical knowledge,
we do no longer provide the specific citations. The main existence result is
provided in Section 2.

\subsection{Preliminaries}

The indicator function $I_{K}:K\rightarrow (-\infty ,\infty ]$ equals zero
at a point of $K$ and $+\infty $ otherwise. Then, $\partial
I_{K}:K\rightarrow 2^{\mathbb{R}}$ is defined by%
\begin{equation*}
\partial I_{K}(\zeta )=\{\chi \in \mathbb{R};\mbox{ }\chi (\zeta -\overline{%
\zeta })\geq 0,\mbox{ }\forall \overline{\zeta }\in K\}.
\end{equation*}%
We recall that $\partial I_{K}(\zeta )=N_{K}(\zeta )\subset \mathbb{R},$ the
normal cone to $K$ at $\zeta .$ In particular, we have%
\begin{equation*}
N_{K}(\zeta )=\{\chi \in \mathbb{R};\mbox{ }\chi \leq 0\mbox{ if }\zeta
=\Gamma ^{\ast }\mbox{ and }\chi =0\mbox{ if }\zeta >\Gamma ^{\ast }\}.
\end{equation*}%
Also, we denote by $\mathcal{K}$ the set 
\begin{equation}
\mathcal{K}=\{z\in L^{\infty }(0,T);\mbox{ }z(t)\geq \Gamma ^{\ast }\mbox{
a.e. }t\in (0,T)\}  \label{L1-5}
\end{equation}%
and by $\mathcal{N}_{\mathcal{K}}(z)\subset (L^{\infty }(0,T))^{\ast }$ the
corresponding normal cone at $z\in \mathcal{K},$ that is%
\begin{equation}
\mathcal{N}_{\mathcal{K}}(z)=\{\eta \in (L^{\infty }(0,T))^{\ast };\mbox{ }%
\eta (z-y)\geq 0,\mbox{ }\forall y\in \mathcal{K}\},  \label{L1-6}
\end{equation}%
where $\eta (z-y)$ is the value of the functional $\eta \in (L^{\infty
}(0,T))^{\ast }$ at $(z-y)\in L^{\infty }(0,T)$ (see \cite{vbp-2012}, p.
242-244). We note that $v\in \mathcal{N}_{\mathcal{K}}(z)\cap L^{1}(0,T)$
iff $v(t)\in N_{K}(z(t)),$ a.e. $t\in (0,T).$

Here, $(L^{\infty }(0,T))^{\ast }$ is the dual of the space $L^{\infty
}(0,T) $ of essentially bounded functions on $(0,T).$ The space $(L^{\infty
}(0,T))^{\ast }$ is a linear subspace of $\mathcal{M}([0,T])=(C[0,T])^{\ast
},$ the space of bounded Radon measures on $[0,T].$ We note that by the
Lebesgue decomposition theorem (see e.g. \cite{Rockafeller-II-71}), every $%
\mu \in (L^{\infty }(0,T))^{\ast }$ can be uniquely written as 
\begin{equation}
\mu =\mu _{a}+\mu _{s},  \label{L1-2}
\end{equation}%
where $\mu _{a}\in L^{1}(0,T)$ and $\mu _{s}$ is a singular measure (that is
there exists a Lebesgue measurable set $S\subset \lbrack 0,T]$ with meas$%
([0,T]\backslash S)=0$ and $\mu _{s}(\varphi )=0$ for all $\varphi \in
L^{\infty }(S)).$ This means that $\mu _{s}$ has the support on a set of
zero measure ($[0,T]\backslash S).$

We denote by $BV([0,T])$ the space of functions $v:[0,T]\rightarrow \mathbb{R%
}$ with bounded variation, that is 
\begin{equation*}
\left\Vert v\right\Vert _{BV([0,T])}=\sup \left\{
\sum_{i=0}^{N-1}\left\vert v(t_{i+1})-v(t_{i})\right\vert ;\mbox{ }%
0=t_{0}<t_{1}<...<t_{N}=T\right\} <\infty .
\end{equation*}%
Every $v\in BV([0,T])$ has a unique decomposition (see e.g., \cite{Ambrosio}%
) 
\begin{equation}
v=v^{a}+v^{s},  \label{L1-3}
\end{equation}%
where $v^{a}\in AC[0,T]$ and $v^{s}\in BV([0,T]).$ Here, $AC[0,T]$ is the
space of absolutely continuous functions on $[0,T]$ and $v^{s}$ is a
singular part (for instance it can be a jump function with bounded variation
or a function with bounded variation with a.e. zero derivative).

\noindent We note that if $v\in BV([0,T]),$ then its distributional
derivative $\frac{dv}{dt}:=\mu $ belongs to $(L^{\infty }(0,T))^{\ast },$
and in virtue of the Lebesgue decomposition, we have the following
representation, as the sum of the absolutely continuous part (in the sense
of measure) and the singular part%
\begin{equation}
\frac{dv}{dt}=\mu _{a}+\mu _{s}=\frac{dv^{a}}{dt}+\frac{dv^{s}}{dt}\in 
\mathcal{D}^{\prime }(0,T),  \label{L1-4}
\end{equation}%
where $\mathcal{D}^{\prime }(0,T)$ is the space of Schwartz distributions on 
$(0,T).$

\section{The main result}

\setcounter{equation}{0}

We shall assume that the following hypotheses hold for the functions
occurring in (\ref{L1}):%
\begin{eqnarray}
U &\in &C([0,\infty )\times \lbrack 0,\infty )),\mbox{ }  \label{L2} \\
x &\rightarrow &U(t,x)\mbox{ is Lipschitz continuous with the Lipschitz
constant }U_{Lip},\mbox{ for all }t\geq 0,  \notag \\
0 &\leq &U(t,x),\mbox{ for all }(t,x)\in \lbrack 0,\infty )\times \lbrack
0,\infty ),\mbox{ }  \notag \\
0 &\leq &U(t,0)\leq U_{0\max },\mbox{ for }t\geq 0,  \notag
\end{eqnarray}%
\begin{eqnarray}
&&\Gamma \in C([0,\infty )\times \lbrack 0,\infty )),\mbox{ }\frac{\partial
\Gamma }{\partial t}\in L^{\infty }((0,\infty )\times (0,\infty )),
\label{L3} \\
&&0\leq \Gamma (t,x)\leq \Gamma _{\max },\mbox{ for all }(t,x)\in \lbrack
0,\infty )\times \lbrack 0,\infty ),\mbox{ \ }  \notag \\
&&\left\vert \frac{\partial \Gamma }{\partial t}(t,x)\right\vert \leq
C_{\Gamma },\mbox{ for all }(t,x)\in (0,\infty )\times \lbrack 0,\infty )%
\mbox{,}  \notag
\end{eqnarray}%
and 
\begin{equation}
\eta _{\ast }:=\inf_{y\geq L_{0}}\left\{ \frac{1}{y}\int_{L_{0}}^{y}\Gamma
(t,\sigma )d\sigma \right\} ,\mbox{ }\eta _{\ast }>\Gamma ^{\ast }>0,\mbox{ }%
L_{0}\geq 0.  \label{L3-1}
\end{equation}%
An example of a function $\Gamma $ complying with these hypotheses is of the
form $\Gamma (t,\sigma )=a(t)-2\gamma _{0}\sigma e^{-\sigma ^{2}},$ where $%
a(t)>\eta _{\ast }$ and $\gamma _{0}>0.$ Then, 
\begin{equation*}
\inf \left\{ \frac{1}{y}\int_{0}^{y}\Gamma (t,\sigma )d\sigma \right\} =\inf
\left\{ a(t)-\gamma _{0}\frac{1-e^{-y^{2}}}{y}\right\} =a(t)=\eta _{\ast
}>\Gamma ^{\ast }.
\end{equation*}%
In particular, for $a$ a positive constant, the graphic of such a function
may describe well the behavior of the cohesion function in the cell growth
model discussed in the introduction.

$\medskip $

\noindent \textbf{Definition 2.1. }The function $L:[0,T]\rightarrow \mathbb{R%
}$ is called a solution to the variational inequality (\ref{L1}) if the
following conditions hold:

$\medskip $

\begin{equation}
L\in BV([0,T]),\mbox{ }L(0)=L_{0},\mbox{ }L=L^{a}+L^{s},  \label{L1-0}
\end{equation}

\begin{eqnarray}
\frac{dL^{a}}{dt}+\mu _{a}(t) &=&U(t,L(t)),\mbox{ a.e. }t\in (0,T),
\label{L1-7} \\
\frac{dL^{s}}{dt}+\mu _{s} &=&0\mbox{, in }\mathcal{D}^{\prime }(0,T), 
\notag
\end{eqnarray}

\begin{eqnarray}
\mu _{a} &\in &L^{1}(0,T),\mbox{ }\mu _{s}\in (L^{\infty }(0,T))^{\ast },
\label{L1-8} \\
\mu _{a}(t) &\in &N_{K}(\Gamma (t,L(t))),\mbox{ a.e. }t\in (0,T),\mbox{ }\mu
_{s}\in \mathcal{N}_{\mathcal{K}}(\Gamma (\cdot ,L(\cdot ))).  \notag
\end{eqnarray}

\medskip

\noindent Here, $L^{a}\in AC[0,T]$ is the absolutely continuous part of $L$
and $L^{s}$ is the singular part, while $\mu _{a}$ and $\mu _{s}$ are the
absolutely continuous and singular parts, respectively of $\mu \in
(L^{\infty }(0,T))^{\ast }$.

\medskip

\noindent \textbf{Theorem 2.2.} \textit{Under the assumptions} (\ref{L2})-(%
\ref{L3-1}), \textit{the variational inequality }(\ref{L1})\textit{\ has at
least one solution, }$L\in BV([0,T]),$ \textit{that is }$L=L^{\alpha
}+L^{s}, $\textit{\ satisfying the equations}%
\begin{equation}
\frac{dL^{a}}{dt}(t)+\mu _{a}(t)=U(t,L(t)),\mbox{ \textit{a.e}. }t\in (0,T),
\label{La}
\end{equation}%
\begin{equation}
\frac{dL^{s}}{dt}+\mu _{s}=0,\mbox{ \ \ \ \ \textit{in} }\mathcal{D}^{\prime
}(0,T),  \label{Ls}
\end{equation}%
\textit{where} 
\begin{equation}
\mu _{a}(t)\in N_{K}(\Gamma (t,L(t)),\mbox{ \textit{a.e}. }t\in (0,T),
\label{miua}
\end{equation}%
\textit{and}%
\begin{equation}
\mu _{s}\in \mathcal{N}_{\mathcal{K}}(\Gamma (\cdot ,L(\cdot ))).
\label{mius}
\end{equation}

\medskip

\noindent \textbf{Proof.} For $\varepsilon >0$ we introduce the Yosida
approximation of $\partial I_{K},$ 
\begin{equation*}
(\partial I_{K})_{\varepsilon }(z)=\frac{1}{\varepsilon }(I-(I+\varepsilon
\partial I_{K})^{-1})z\mbox{, }
\end{equation*}%
and denote by $P_{K}z$ the projection of $z\in \mathbb{R}$ on $K,$ given by 
\begin{equation*}
P_{K}z=\left\{ 
\begin{array}{l}
z,\mbox{ \ \ }z>\Gamma ^{\ast } \\ 
\Gamma ^{\ast },\mbox{ }z\leq \Gamma ^{\ast }.%
\end{array}%
\right.
\end{equation*}%
We recall that 
\begin{equation*}
(I+\varepsilon \partial I_{K})^{-1}z=P_{K}z,\mbox{ for all }z\in \mathbb{R}%
\mbox{.}
\end{equation*}%
This implies that 
\begin{equation*}
(\partial I_{K})_{\varepsilon }(z)=\left\{ 
\begin{array}{l}
0,\mbox{ \ \ \ \ \ \ \ \ \ \ \ }z>\Gamma ^{\ast } \\ 
\frac{1}{\varepsilon }(z-\Gamma ^{\ast }),\mbox{ }z\leq \Gamma ^{\ast }%
\end{array}%
\right. =-\frac{1}{\varepsilon }(z-\Gamma ^{\ast })^{-},\mbox{ }\forall z\in 
\mathbb{R}\mbox{.}
\end{equation*}%
Here, $(\cdot )^{-}$ represents the negative part.

Let us consider the approximating problem 
\begin{eqnarray}
\frac{dL_{\varepsilon }}{dt}(t)+(\partial I_{K})_{\varepsilon }(\Gamma
(t,L_{\varepsilon }(t))) &=&U(t,L_{\varepsilon }(t)),\mbox{ a.e. }t\in (0,T),
\label{L4} \\
L_{\varepsilon }(0) &=&L_{0}\geq 0.  \notag
\end{eqnarray}%
It is obvious that (\ref{L4}) has a unique solution $L_{\varepsilon }\in
C^{1}[0,T]$ satisfying 
\begin{equation*}
L_{\varepsilon }(t)=L_{0}+\int_{0}^{t}(U(s,L_{\varepsilon }(s))-\mu
_{\varepsilon }(s))ds,\mbox{ for all }t\in \lbrack 0,T],
\end{equation*}%
where, 
\begin{equation}
\mu _{\varepsilon }(t):=(\partial I_{K})_{\varepsilon }(z_{\varepsilon }(t)),%
\mbox{ \ }z_{\varepsilon }(t):=\Gamma (t,L_{\varepsilon }(t)).  \label{L4-1}
\end{equation}%
Moreover, since $U-\mu _{\varepsilon }$ is positive it follows that $%
L_{\varepsilon }(t)\geq L_{0}$, for all $t\in \lbrack 0,T].$

\noindent To continue the proof we need some estimates. We multiply (\ref{L4}%
) by $(\Gamma (s,L_{\varepsilon }(s))-\alpha ),$ where $\alpha $ is a
positive constant which will be specified a little later, and integrate on $%
(0,t)$ to obtain 
\begin{eqnarray}
&&\int_{0}^{t}\frac{dL_{\varepsilon }}{ds}(s)\Gamma (s,L_{\varepsilon
}(s))ds-\int_{0}^{t}\alpha \frac{dL_{\varepsilon }}{ds}(s)ds+\int_{0}^{t}\mu
_{\varepsilon }(s)(z_{\varepsilon }(s)-\alpha )ds  \label{L7} \\
&=&\int_{0}^{t}U(s,L_{\varepsilon }(s))(\Gamma (s,L_{\varepsilon
}(s))-\alpha )ds.  \notag
\end{eqnarray}%
We denote by $j$ the potential of the function $\sigma \rightarrow \Gamma
(t,\sigma ),$ that is 
\begin{equation}
j(t,v):=\int_{L_{0}}^{v}\Gamma (t,\sigma )d\sigma ,\mbox{ }v\in \mathbb{R},%
\mbox{ }v\geq L_{0},\mbox{ for all }t\in \lbrack 0,T],  \label{L8}
\end{equation}%
and note that%
\begin{equation*}
\frac{d}{dt}\int_{L_{0}}^{L_{\varepsilon }(t)}\Gamma (t,\sigma )d\sigma =%
\frac{dL_{\varepsilon }}{dt}(t)\Gamma (t,L_{\varepsilon
}(t))+\int_{L_{0}}^{L_{\varepsilon }(t)}\frac{\partial \Gamma }{\partial t}%
(t,\sigma )d\sigma ,\mbox{ }t\in \lbrack 0,T].
\end{equation*}%
This yields%
\begin{equation*}
\int_{0}^{t}\frac{dL_{\varepsilon }}{ds}(s)\Gamma (s,L_{\varepsilon
}(s))ds=\int_{L_{0}}^{L_{\varepsilon }(t)}\Gamma (t,\sigma )d\sigma
-\int_{0}^{t}\int_{L_{0}}^{L_{\varepsilon }(s)}\frac{\partial \Gamma }{%
\partial s}(s,\sigma )d\sigma ds.
\end{equation*}%
Replacing the left-hand side term of the previous equality in (\ref{L7}) we
get%
\begin{eqnarray}
&&j(t,L_{\varepsilon }(t))-\int_{0}^{t}\int_{L_{0}}^{L_{\varepsilon }(s)}%
\frac{\partial \Gamma }{\partial s}(s,\sigma )d\sigma ds-L_{\varepsilon
}(t)\alpha +L_{0}\alpha  \notag \\
&&+\int_{0}^{t}\mu _{\varepsilon }(s)(z_{\varepsilon }(s)-\alpha
)ds=\int_{0}^{t}U(s,L_{\varepsilon }(s))(z_{\varepsilon }(s)-\alpha )ds.
\label{L8-1}
\end{eqnarray}%
By (\ref{L3-1}) we have 
\begin{equation}
\eta _{\ast }=\inf_{v\geq L_{0}}\left\{ \frac{1}{v}\int_{L_{0}}^{v}\Gamma
(t,\sigma )d\sigma \right\} \leq \frac{1}{v}\int_{L_{0}}^{v}\Gamma (t,\sigma
)d\sigma ,\mbox{ for }v\geq L_{0},  \label{L8-2}
\end{equation}%
and so for $v=L_{\varepsilon }(t)$ we get 
\begin{equation*}
\eta _{\ast }L_{\varepsilon }(t)\leq \int_{L_{0}}^{L_{\varepsilon
}(t)}\Gamma (t,\sigma )d\sigma ,\mbox{ }\forall t\in \lbrack 0,T].
\end{equation*}%
Using (\ref{L2})-(\ref{L3-1}) and (\ref{L8-2}) in (\ref{L8-1}) we obtain%
\begin{eqnarray*}
&&(\eta _{\ast }-\alpha )L_{\varepsilon }(t)+\int_{0}^{t}\mu _{\varepsilon
}(s)(z_{\varepsilon }(s)-\alpha )ds \\
&\leq &\int_{0}^{t}\int_{L_{0}}^{L_{\varepsilon }(s)}\left\vert \frac{%
\partial \Gamma }{\partial s}(s,\sigma )\right\vert d\sigma
ds+\int_{0}^{t}\left\vert U(s,L_{\varepsilon }(s))\right\vert \left\vert
z_{\varepsilon }(s)-\alpha \right\vert ds \\
&\leq &C_{\Gamma }\int_{0}^{t}L_{\varepsilon }(s)ds-C_{\Gamma
}L_{0}T+(\Gamma _{\max }+\alpha )\int_{0}^{t}\left\vert
U_{0}(s)+U_{Lip}L_{\varepsilon }(s))\right\vert ds,
\end{eqnarray*}%
where $U_{0}(t)=U(t,0)$ and $U_{Lip}$ is the Lipschitz constant of $U.$
Recall that $z_{\varepsilon }(s)=\Gamma (s,L_{\varepsilon }(s))\leq \Gamma
_{\max }.$ This yields%
\begin{equation}
(\eta _{\ast }-\alpha )L_{\varepsilon }(t)+\int_{0}^{t}(\partial
I_{K})_{\varepsilon }(\Gamma (s,L_{\varepsilon }(s)))(\Gamma
(s,L_{\varepsilon }(s))-\alpha )ds\leq C_{1}+C_{2}\int_{0}^{t}L_{\varepsilon
}(s)ds,  \label{L9}
\end{equation}%
where $C_{1}(T)=(\Gamma _{\max }+\alpha )U_{0\max }T,$ $C_{2}=(\Gamma _{\max
}+\alpha )U_{Lip}+C_{\Gamma }.$

\noindent Now, we can choose $\alpha <\eta _{\ast }$ and so, for 
\begin{equation}
\eta _{\ast }>\alpha >\Gamma ^{\ast }  \label{L9-1}
\end{equation}%
the first term on the left-hand side is positive. Applying the Gronwall
lemma we obtain 
\begin{equation}
L_{\varepsilon }(t)\leq \frac{C_{1}}{\eta _{\ast }-\alpha }e^{\frac{C_{2}}{%
\eta _{\ast }-\alpha }t}\mbox{ for all }t\in \lbrack 0,T].  \label{L10}
\end{equation}%
Also, it follows that 
\begin{equation}
0\leq \int_{0}^{t}(\partial I_{K})_{\varepsilon }(\Gamma (s,L_{\varepsilon
}(s)))(\Gamma (s,L_{\varepsilon }(s))-\alpha )ds\leq C_{1}e^{\frac{C_{2}}{%
\eta _{\ast }-\alpha }t},\mbox{ for all }t\in \lbrack 0,T].  \label{L11}
\end{equation}%
By the definition of the subdifferential, we are entitled to write 
\begin{equation}
\mu _{\varepsilon }(t)(z_{\varepsilon }(t)-\alpha -\rho \theta )\geq 0,\mbox{
a.e. }t\in (0,T),  \label{L5}
\end{equation}%
where $\alpha ,$ $\rho $ and $\theta $ are such that $\alpha +\rho \theta
\geq \Gamma ^{\ast },$ with $\alpha >\Gamma ^{\ast },$ $\theta <0,$ $%
\left\vert \theta \right\vert =1,$ $0<\rho \leq \alpha -\Gamma ^{\ast }.$
Now, if $1_{A}$ is the characteristic function of the set $A,$ we set 
\begin{equation*}
\theta (t):=\frac{\mu _{\varepsilon }(t)}{\left\vert \mu _{\varepsilon
}(t)\right\vert }1_{\{t;\mbox{ }\mu _{\varepsilon }(t)\neq 0\}}
\end{equation*}%
and by integrating (\ref{L5}), we get%
\begin{equation}
\rho \int_{0}^{t}\left\vert \mu _{\varepsilon }(s)\right\vert ds\leq
\int_{0}^{t}\mu _{\varepsilon }(s)(z_{\varepsilon }(s)-\alpha )ds.
\label{L6}
\end{equation}%
By (\ref{L11}) we deduce that%
\begin{equation}
\int_{0}^{t}\left\vert (\partial I_{K})_{\varepsilon }(\Gamma
(s,L_{\varepsilon }(s)))\right\vert ds\leq \frac{C_{1}}{\rho }e^{\frac{C_{2}%
}{\eta _{\ast }-\alpha }t},\mbox{ for all }t\in \lbrack 0,T],  \label{L12}
\end{equation}%
while by (\ref{L4}) we obtain 
\begin{eqnarray}
\int_{0}^{t}\left\vert \frac{dL_{\varepsilon }}{ds}(s)\right\vert ds &\leq
&\int_{0}^{t}\left\vert U(s,L_{\varepsilon }(s))\right\vert
ds+\int_{0}^{t}\left\vert (\partial I_{K})_{\varepsilon }(\Gamma
(s,L_{\varepsilon }(s)))\right\vert ds  \label{L13} \\
&\leq &\int_{0}^{t}(U_{0\max }+U_{Lip}\left\vert L_{\varepsilon
}(s)\right\vert )ds+\frac{C_{1}}{\rho }e^{\frac{C_{2}}{\eta _{\ast }-\alpha }%
t}  \notag \\
&\leq &U_{0\max }T+\frac{C_{1}}{\rho }e^{\frac{C_{2}}{\eta _{\ast }-\alpha }%
T}+U_{Lip}\frac{C_{1}}{C_{2}}\left( e^{\frac{C_{2}}{\eta _{\ast }-\alpha }%
t}-1\right) ,\mbox{ for all }t\in \lbrack 0,T].  \notag
\end{eqnarray}%
Writing (\ref{L4}) as 
\begin{equation*}
\frac{dL_{\varepsilon }}{dt}(t)-\frac{1}{\varepsilon }(\Gamma
(t,L_{\varepsilon }(t))-\Gamma ^{\ast })^{-}=U(t,L_{\varepsilon }(t))
\end{equation*}%
and multiplying by $(\Gamma (t,L_{\varepsilon }(t))-\Gamma ^{\ast })$ we get 
\begin{equation*}
\frac{1}{\varepsilon }\left\Vert (\Gamma (\cdot ,L_{\varepsilon }(\cdot
))-\Gamma ^{\ast })^{-}\right\Vert _{L^{2}(0,T)}^{2}=\int_{0}^{t}\left(
U(s,L_{\varepsilon }(s))-\frac{dL_{\varepsilon }}{ds}(s)\right) (\Gamma
(s,L_{\varepsilon }(s))-\Gamma ^{\ast })ds\leq C_{T}.
\end{equation*}%
This immediately yields 
\begin{equation}
(\Gamma (t,L_{\varepsilon }(t))-\Gamma ^{\ast })^{-}\leq C_{T}\varepsilon ,%
\mbox{ for all }t\in \lbrack 0,T].  \label{L13-1}
\end{equation}%
By $C_{T}$ we denote several constants, which can differ from line to line.
They depend on the data and $T,$ but are independent of $\varepsilon .$

We conclude that $\{L_{\varepsilon }\}_{\varepsilon }$ is bounded in $C[0,T]$%
, $\left\{ \frac{dL_{\varepsilon }}{dt}\right\} _{\varepsilon }$ and $%
\left\{ \mu _{\varepsilon }=(\partial I_{K})_{\varepsilon }(\Gamma (\cdot
,L_{\varepsilon }(\cdot ))\right\} _{\varepsilon }$ are bounded in $%
L^{1}(0,T).$

By (\ref{L13}) it follows that $\left\Vert L_{\varepsilon }\right\Vert
_{BV([0,T])}\leq C_{T}$ for all $\varepsilon >0,$ and so, by Helly's theorem
(see e.g., \cite{vbp-2012}, p. 47), we have $L\in BV([0,T])$ and%
\begin{equation}
L_{\varepsilon }(t)\rightarrow L(t),\mbox{ for all }t\in \lbrack 0,T].
\label{L14}
\end{equation}%
In particular, $L_{\varepsilon }(0)\rightarrow L(0)=L_{0}.$ Moreover, by
Egorov's theorem, for each $\delta >0,$ there exists a set $\Omega _{\delta
}\subset \lbrack 0,T],$ such that meas$(\Omega _{\delta })<\delta $ and 
\begin{equation}
L_{\varepsilon }\rightarrow L\mbox{ uniformly on }\Omega _{\delta },\mbox{
as }\varepsilon \rightarrow 0.  \label{L14-0}
\end{equation}%
Then, the sequences $\left\{ \frac{dL_{\varepsilon }}{dt}\right\}
_{\varepsilon }$ and $\left\{ \mu _{\varepsilon }\right\} _{\varepsilon }$
are weak-* compact in $(L^{\infty }(0,T))^{\ast },$ as specified in the
proof of Corollary 2B in \cite{Rockafeller-II-71}. We stress that this is
not directly implied by the Alaoglu theorem, but can be deduced by the
following argument. Let us consider the linear operator $\Phi
:C[0,T]\rightarrow L^{\infty }(Q),$ $\Phi v=\widetilde{\Phi },$ which maps a
continuous function into the corresponding class of equivalence $\widetilde{%
\Phi }$ (of all functions a.e. equal). Its adjoint $\Phi ^{\ast }:(L^{\infty
}(Q))^{\prime }\rightarrow \mathcal{M}([0,T])$ is defined by $(\Phi ^{\ast
}\mu )(v):=\mu (\Phi v)$ for any $v\in C[0,T].$ If $\{\mu _{n}\}_{n}$ is
bounded in $(L^{\infty }(Q))^{\prime }$ and also in $\mathcal{M}([0,T]),$
then $\{\Phi ^{\ast }\mu _{n}\}_{n}$ is bounded in $\mathcal{M}([0,T])$
which is the dual of the separable space $C[0,T]$ and so by the Alaoglu
theorem $\{\Phi ^{\ast }\mu _{n}\}_{n}$ is weak-*\ sequentially compact in $%
\mathcal{M}([0,T]).$ Also $\{\mu _{n}\}_{n}$ is weak-*\ sequentially compact
in $\mathcal{M}([0,T]).$ Passing to the limit in $\mu _{n}(\Phi v)=(\Phi
^{\ast }\mu _{n})(v)$ we get $\mu (\Phi v)=(\Phi ^{\ast }\mu )(v):=$ for any 
$\widetilde{\Phi }\in L^{\infty }(Q)$ which is of the form $\Phi v$ with $%
v\in C[0,T].$ Then, $\mu $ can be extended by the Hahn-Banach theorem to all 
$L^{\infty }(0,T)$ and so we conclude that $\{\mu _{n}\}_{n}$ is weak-star\
sequentially compact in $(L^{\infty }(Q))^{\prime }.$ Therefore, one can
extract a subsequence such that%
\begin{equation}
\frac{dL_{\varepsilon }}{dt}\rightarrow \frac{dL}{dt}\mbox{ weak-* in }%
(L^{\infty }(0,T))^{\ast }\subset \mathcal{M}([0,T]),  \label{L14-1}
\end{equation}%
\begin{equation}
\mu _{\varepsilon }\rightarrow \mu \mbox{ weak-* in }(L^{\infty
}(0,T))^{\ast }\subset \mathcal{M}([0,T]).  \label{L14-2}
\end{equation}%
Since $\Gamma (t,x)$ and $U(t,x)$ are continuous with respect to $x$ it
follows that 
\begin{equation}
\Gamma (t,L_{\varepsilon }(t))\rightarrow \Gamma (t,L(t)),\mbox{ for all }%
t\in \lbrack 0,T],  \label{L15}
\end{equation}%
\begin{equation}
U(t,L_{\varepsilon }(t))\rightarrow U(t,L(t)),\mbox{ for all }t\in \lbrack
0,T],  \label{L16}
\end{equation}%
and by (\ref{L13-1}) we get 
\begin{equation}
\Gamma (t,L(t))\geq \Gamma ^{\ast },\mbox{ for all }t\in \lbrack 0,T].
\label{L16-0}
\end{equation}%
Moreover, by (\ref{L4}) we have at limit 
\begin{equation}
\frac{dL}{dt}+\mu =U(\cdot ,L(\cdot ))\mbox{ in }\mathcal{D}^{\prime }(0,T).
\label{L17}
\end{equation}%
Further, by writing 
\begin{equation*}
\int_{\Omega _{\delta }}\mu _{\varepsilon }(t)(\Gamma (t,L_{\varepsilon
}(t))-v(t))dt=\int_{0}^{T}\mu _{\varepsilon }(t)\left( 1_{\Omega _{\delta
}}(t)(\Gamma (t,L_{\varepsilon }(t))-v(t))\right) dt\geq 0\mbox{,}
\end{equation*}%
for all $v\in L^{\infty }(0,T),$ $v(t)\geq \Gamma ^{\ast }$ a.e., we obtain
at limit $\mu (1_{\Omega _{\delta }}(\Gamma (\cdot ,L(\cdot ))-v))\geq 0$
for all $v\in \mathcal{K}$. Since $\delta $ is positive arbitrary, we get as 
$\delta \rightarrow 0$ 
\begin{equation}
\mu (\Gamma (\cdot ,L(\cdot ))-v))\geq 0,\mbox{ for all }v\in \mathcal{K},
\label{L18}
\end{equation}%
and so, $\mu \in \mathcal{N}_{\mathcal{K}}(\Gamma (\cdot ,L(\cdot ))).$

Since $\mu \in (L^{\infty }(0,T))^{\ast }\subset \mathcal{M}([0,T]),$ $\mu $
can be written as $\mu =\mu _{a}+\mu _{s},$ where $\mu _{a}$ is the
absolutely continuous part (in the sense of measure) and $\mu _{s}$ is the
singular part and so (\ref{L17}) implies that%
\begin{equation*}
\frac{dL^{a}}{dt}(t)+\mu _{a}(t)=U(t,L(t)),\mbox{ a.e. }t\in (0,T),
\end{equation*}%
\begin{equation*}
\frac{dL^{s}}{dt}+\mu _{s}=0,\mbox{ \ \ \ \ \ in }\mathcal{D}^{\prime }(0,T),
\end{equation*}%
where 
\begin{equation*}
\mu _{a}(t)\in N_{K}(\Gamma (t,L(t)),\mbox{ a.e. }t\in (0,T),\mbox{ }\mu
_{s}\in \mathcal{N}_{\mathcal{K}}(\Gamma (\cdot ,L(\cdot ))),
\end{equation*}%
as claimed. This means that 
\begin{eqnarray*}
\mu _{s}(\varphi ) &=&0\mbox{ if }\varphi \in \overset{\circ }{\mathcal{K}}
\\
\mu _{s}(\varphi ) &\leq &0\mbox{ if }\varphi \in \partial \mathcal{K}\mbox{.%
}
\end{eqnarray*}%
Recalling that 
\begin{equation*}
\mbox{supp }\mu _{s}=\{\Sigma \subset \lbrack 0,T];\mbox{ }\mu _{s}\neq 0%
\mbox{ on }\Sigma ,\mbox{ i.e., }\mu _{s}(\varphi )\neq 0,\mbox{ for }%
\varphi \in \Sigma \}
\end{equation*}%
it follows that supp $\mu _{s}\subset \partial \mathcal{K}=\{z\in L^{\infty
}(0,T);$ $z(t)=\Gamma ^{\ast }\}$ of $\mathcal{K}.$ Thus, we have in fact 
\begin{equation*}
\mbox{supp }\mu _{s}=\{t\in \lbrack 0,T];\mbox{ }\Gamma (t,L(t))=\Gamma
^{\ast }\}.
\end{equation*}%
Therefore, $L$ has an absolutely continuous part $L^{a}$ and a $BV$ part $%
L^{s},$ where $\Gamma (t,L(t))=\Gamma ^{\ast }.$ In particular, we note that 
$L^{s}$ can be represented as jump functions at $t$, e.g. $L^{s}(t)=\alpha
_{i}$ on $[t_{i},t_{i+1}),$ meaning that these jump points are those at
which $\Gamma (t_{i},L(t_{i}))=\Gamma ^{\ast }.$ This completes the proof of
the solution existence.\hfill $\square $

$\medskip $

\noindent \textbf{Remark 2.3. }It should be noted that since $L\rightarrow
\partial I_{K}(\Gamma (\cdot ,L))$ is not monotone, the uniqueness remains
open.

\end{document}